\documentclass{article}

\usepackage{amssymb}
\usepackage{eucal}
\usepackage{amsmath}
\usepackage{amsthm}

\numberwithin{equation}{section}
\newtheorem{thm}{\bf Theorem}[section]

\theoremstyle{remark}

\title{Hamilton-Poisson formulation for the rotational motion of a rigid body in the presence of an axisymmetric force field and a gyroscopic torque}
\author{Petre Birtea\footnote{Corresponding author - E-mail: birtea@math.uvt.ro; Phone: +40 726 424147; Fax: +40 256 592316}, Ioan Ca\c su, Dan Com\u anescu\\
{\small Department of Mathematics, West University of Timi\c soara}\\
{\small Bd. V. P\^ arvan, No 4, 300223 Timi\c soara, Rom\^ ania}\\
{\small E-mail addresses: birtea@math.uvt.ro; casu@math.uvt.ro; comanescu@math.uvt.ro}}
\date{}

\begin{document}

\maketitle



\begin{abstract}
We give sufficient conditions for the rigid body in the presence of an axisymmetric force field and a gyroscopic torque to admit a Hamilton-Poisson formulation. Even if by adding a gyroscopic torque we initially lose one of the conserved Casimirs, we recover  another conservation law as a Casimir function for a modified Poisson structure. We apply this frame to several well known results in the literature.
\end{abstract}

\section{Introduction}

The generalized Euler-Poisson system, which is a dynamical system for the rotational motion of a rigid body in the presence of an axisymmetric force field, admits a Hamilton-Poisson formulation, where the bracket is the "$-$" Kirillov-Kostant-Souriau bracket on the dual Lie algebra $e(3)^*$. The Hamiltonian function is of the type kinetic energy plus potential energy. The Hamiltonian function and the two Casimir functions for the K-K-S Poisson structure are conservation laws for the dynamic. 

Adding a certain type of gyroscopic torque the Hamiltonian and one of the Casimirs remain conserved along the solutions of the new dynamical system. In this paper we give a technique for finding a third conservation law in the presence of gyroscopic torque. 

We modify the K-K-S structure using the elements that generate the gyroscopic torque and we find necessary and sufficient conditions for this structure to be a Poisson structure. This leads us to what we call the {\bf Jacobi condition equation}. Once this new Poisson structure is established, we give necessary and sufficient conditions for a function to be Casimir function. We call these conditions the {\bf Casimir condition equation}. 

Using this framework, we study several general cases of gyroscopic torques, that include some well known situations found in the literature. We encompass the gyrostatic torques, affine gyroscopic torques, nonlinear torques studied by Yehia \cite{Yeh3}, \cite{Yeh4}, gyroscopic torques along one axis of inertia studied by Borisov and Mamaev \cite{Bor}.

It would be interesting to relate the methods presented in this paper with the methods of constructing and analyzing gyroscopic torques given in Chang and Marsden \cite{Chang1} and in Chang et al. \cite{Chang2}. The understanding of the link between the two methods could give better results on studying the dynamics in the presence of gyroscopic torque.

\section{Mechanical aspects}

The equation for the rotation of a rigid body about a point $O$ is
\begin{equation*}
\mathbb{I}\dot{\vec{\omega}}+\vec{\omega}\times\mathbb{I}\vec{\omega}=\vec{\mathcal{M}}_O
\end{equation*}
in the cases:

\noindent (i) the point $O$ is fixed with respect to an inertial frame;

\noindent (ii) the point $O$ coincides with the center of mass.\\
In the above equation $\vec{\omega}$ is the angular velocity, $\mathbb{I}$ is the moment of inertia tensor at the point $O$ and
$\vec{\mathcal{M}}_O$ is the applied torque about the point $O$.
For the following considerations we use the body frame $\mathcal{R}$ with the origin $O$ and its axes are principal axes of inertia of the
rigid body at the point $O$. The matrix of the moment of inertia tensor in this body frame has the form
\begin{equation*}
\mathbb{I}_{\mathcal{R}}=\hbox{diag}(A,B,C),
\end{equation*}
where $A$, $B$ and $C$ are principal moments of inertia.

Our rigid body is situated in an axisymmetric force field. We denote by $\vec{\gamma}$ a vector which generates the axis of symmetry and we
have the Poisson's differential equation
\begin{equation*}
\dot{\vec{\gamma}}+\vec{\omega}\times\vec{\gamma}=0,
\end{equation*}
with the time derivative taken with respect to the body frame $\mathcal{R}$.

The rigid body is acted, also, by a gyroscopic torque. We study the case in which the applied torque $\vec{\mathcal{M}}_O$ has the form
\begin{equation*}
\vec{\mathcal{M}}_O=\vec{\mathcal{M}}_a+\vec{\mathcal{M}}_g,
\end{equation*}
where $\vec{\mathcal{M}}_a$ is the torque generated by the axisymmetric force field and $\vec{\mathcal{M}}_g$ is the gyroscopic torque. The
torque $\vec{\mathcal{M}}_a$ depends on $\vec{\gamma}$ and in a sufficiently extensive class of cases it is represented in the form
\begin{equation*}
\vec{\mathcal{M}}_a=\vec{\gamma}\times \nabla_{\gamma}U,
\end{equation*}
where $U$ is a smooth function which depends on $\vec{\gamma}$.
The most general form of the gyroscopic torque considered in this paper is
\begin{equation}\label{torque}
\vec{\mathcal{M}}_g=-\vec{\omega}\times \vec{\mu},
\end{equation}
where $\vec{\mu}$ is a vectorial function which depends on $\vec{\omega}$ and $\vec{\gamma}$.

For mathematical reason we use the angular momentum vector $\vec{M}$ to describe our problem. Its connection with the angular velocity is given by
\begin{equation*}
\vec{M}=\mathbb{I}\vec{\omega}.
\end{equation*}

According with the above considerations the system which describes the rotational motion of a rigid body situated in an axisymmetric force field
and in the presence of a gyroscopic torque has the form
\begin{equation}\label{sistem general}
\left\{%
\begin{array}{ll}
\dot{\vec{M}}=-\mathbb{I}^{-1}\vec{M}\times (\vec{M}+\vec{\mu}(\vec{M},\vec{\gamma}))+\vec{\gamma}\times\nabla_{\gamma}U(\vec{\gamma}) \\
\dot{\vec{\gamma}}=\vec{\gamma}\times\mathbb{I}^{-1}\vec{M}.\end{array}%
\right.
\end{equation}

In the absence of the gyroscopic torque, i.e. $\vec{\mu}=\vec{0}$, the system \eqref{sistem general} becomes
\begin{equation}
\label{redus}
\left\{%
\begin{array}{ll}
\dot{\vec{M}}=-\mathbb{I}^{-1}\vec{M}\times \vec{M}+\vec{\gamma}\times\nabla_{\gamma}U(\vec{\gamma}) \\
\dot{\vec{\gamma}}=\vec{\gamma}\times\mathbb{I}^{-1}\vec{M} \\\end{array}%
\right.
\end{equation}
and we obtain what is known as the generalized Euler-Poisson system. For $\vec{\mu}=\vec{0}$ and $U=mg\vec{\xi}\cdot\vec{\gamma}$ one obtains the classical Euler-Poisson problem which describes (on
the leaf $||\vec{\gamma}||=1$) the rotational motion of the rigid body under gravity. The constant $m$ is the mass of the rigid body, $g$ is the
acceleration of gravity and $\vec{\xi}$ is the vector of position of the mass center.

This system has the following Hamilton-Poisson structure.
The Poisson bracket is the "$-$" Kirillov-Kostant-Souriau bracket on the dual Lie algebra $e(3)^*$. In the coordinates $(M_1,M_2,M_3,\gamma_1,\gamma_2,\gamma_3)$ it has the expression
$$\Pi=\left[\begin{array}{cccccc}
0&-M_3&M_2&0&-\gamma_3&\gamma_2\\
M_3&0&-M_1&\gamma_3&0&-\gamma_1\\
-M_2&M_1&0&-\gamma_2&\gamma_1&0\\
0&-\gamma_3&\gamma_2&0&0&0\\
\gamma_3&0&-\gamma_1&0&0&0\\
-\gamma_2&\gamma_1&0&0&0&0
\end{array}\right].$$ 
The Casimir functions are 
$$C_1(\vec{M},\vec{\gamma})=\frac{1}{2}\vec{\gamma}^2$$
and
$$C_2(\vec{M},\vec{\gamma})=\vec{M}\cdot \vec{\gamma}.$$
The Hamiltonian function is
\begin{equation}
\label{Hamiltonian}
H(\vec{M},\vec{\gamma})=\frac{1}{2}\vec{M}\cdot \mathbb{I}^{-1}\vec{M}+U(\vec{\gamma}).
\end{equation}

In the presence of a gyroscopic torque of the form \eqref{torque} or equivalently $\vec{\mathcal{M}}_g=-\mathbb{I}^{-1}\vec{M}\times \vec{\mu}$ the above functions $H$ and $C_1$ remain constants of motion along the solutions of \eqref{sistem general}.
The time derivative of $C_2$ along the solutions of \eqref{sistem general} is 
$$\frac{d}{dt}C_2=\vec{\gamma}\cdot (\vec{\mu} \times \mathbb{I}^{-1}\vec{M}).$$

Unless $\vec{\mu}$ is a linear combination of $\vec{\gamma}$ and $\mathbb{I}^{-1}\vec{M}=\vec{\omega}$, the function $C_2$ is not a constant of motion for the system \eqref{sistem general}.

\section{Hamilton-Poisson formulation in the presence of gyroscopic torque}

As we have seen in the previous section, when we add a gyroscopic torque to the Hamilton-Poisson system \eqref{redus} we lose one of the Casimirs as a constant of motion for the system \eqref{sistem general}.

We will recover a third constant of motion as a Casimir function by modifying the K-K-S Poisson structure. Our modification of the Poisson structure will be such that the system \eqref{sistem general} becomes a Hamilton-Poisson system.
We consider the following skew-symmetric $6\times 6$ matrix
$$\Pi_{\vec{\mu}}=\left[\begin{array}{cccccc}
0&-M_3-\mu_3&M_2+\mu_2&0&-\gamma_3&\gamma_2\\
M_3+\mu_3&0&-M_1-\mu_1&\gamma_3&0&-\gamma_1\\
-M_2-\mu_2&M_1+\mu_1&0&-\gamma_2&\gamma_1&0\\
0&-\gamma_3&\gamma_2&0&0&0\\
\gamma_3&0&-\gamma_1&0&0&0\\
-\gamma_2&\gamma_1&0&0&0&0
\end{array}\right],$$ 
where the components $\mu_1,\mu_2,\mu_3$ of the vector field $\vec{\mu}$ are smooth functions of the variables $\vec{M}$ and $\vec{\gamma}$. The matrix $\Pi_{\vec{\mu}}$ generates a Poisson structure if and only if the Jacobi identity is verified, i.e.
$$\Pi_{\vec{\mu}}^{li}\partial_l\Pi_{\vec{\mu}}^{jk}+\Pi_{\vec{\mu}}^{lj}\partial_l\Pi_{\vec{\mu}}^{ki}+\Pi_{\vec{\mu}}^{lk}\partial_l\Pi_{\vec{\mu}}^{ij}=0,$$
for all $i,j,k=\overline{1,6}$.
The above system of partial differential equations reduces, in our case, to the following system of ten PDEs for the functions $\mu_1,\mu_2,\mu_3$:
\begin{equation}
\label{sistem Jacobi}
\left\{\begin{array}{l}
\vec{\gamma}\times \nabla_{\vec{M}}\mu_k=\vec{0},~~~k=1,2,3\\
\vec{\gamma}\cdot \hbox{curl}_{\vec{\gamma}}\vec{\mu}+\vec{M}\cdot \hbox{curl}_{\vec{M}}\vec{\mu}+\vec{\mu}\cdot \hbox{curl}_{\vec{M}}\vec{\mu}=0.
\end{array}\right.
\end{equation}

Solving independently the systems of equations $\vec{\gamma}\times \nabla_{\vec{M}}\mu_k=\vec{0}$ for $k=1,2,3$ we obtain $\mu_k$ as functions of four variables $\gamma_1,\gamma_2,\gamma_3$ and $s:=\vec{M}\cdot\vec{\gamma}$,
$$\mu_k=\mu_k(\vec{\gamma},\vec{M}\cdot \vec{\gamma}).$$
Substituting in the last equation of the system \eqref{sistem Jacobi} the $\mu_k$s found above we obtain the {\bf Jacobi condition equation for $\vec{\mu}$}:
\begin{equation}
\label{Jacobi condition}
\vec{\gamma}\cdot \hbox{curl}_{\vec{\gamma}}\vec{\mu}+\vec{\mu}\cdot (\vec{\gamma}\times \partial_4\vec{\mu})=0.
\end{equation}
The partial derivative $\partial_4\vec{\mu}$ is the derivative with respect to the variable $s=\vec{M}\cdot \vec{\gamma}$.
If $\vec{\mu}$ verifies \eqref{Jacobi condition}, then $\Pi_{\vec{\mu}}$ generates a Poisson bracket on $\mathbb{R}^6$.

A smooth function $C(\vec{M},\vec{\gamma})$ is a Casimir function for the Poisson structure $\Pi_{\vec{\mu}}$ if and only if it satisfies the following system of six partial differential equations:
\begin{equation}\label{sistem Casimir}
\left\{\begin{array}{l}
\vec{\gamma}\times \nabla_{\vec{M}}C=\vec{0}\\
(\vec{M}+\vec{\mu})\times \nabla_{\vec{M}}C+\vec{\gamma}\times \nabla_{\vec{\gamma}}C=\vec{0}.
\end{array}\right.
\end{equation}

Solving the system $\vec{\gamma}\times \nabla_{\vec{M}}C=\vec{0}$ we obtain $C$ as a function of four variables $\gamma_1,\gamma_2,\gamma_3,s=\vec{M}\cdot \vec{\gamma}$, i.e. $C=C(\vec{\gamma},\vec{M}\cdot \vec{\gamma})$.
Substituting $C$ found above in the system $(\vec{M}+\vec{\mu})\times \nabla_{\vec{M}}C+\vec{\gamma}\times \nabla_{\vec{\gamma}}C=\vec{0}$ we obtain the {\bf Casimir condition equation}:
\begin{equation}
\label{Casimir condition}
\vec{\gamma}\times ((\partial_4C) \vec{\mu}-\nabla_{\vec{\gamma}}C)=\vec{0}.
\end{equation}
We note that $C_1(\vec{M},\vec{\gamma})=\displaystyle \frac{1}{2}\vec{\gamma}^2$ remains a Casimir for $\Pi_{\vec{\mu}}$ as it verifies the {\bf Casimir condition equation} \eqref{Casimir condition}.
Next, we summarize the considerations we have made in this section.
\begin{thm}
The following results hold:
\begin{itemize}
\item[(i)] $\Pi_{\vec{\mu}}$ generates a Poisson bracket $\{\cdot,\cdot\}_{\vec{\mu}}$ if and only if the vectorial function $\vec{\mu}$ depends only on the four variables $\gamma_1,\gamma_2,\gamma_3,\vec{M}\cdot \vec{\gamma}$ and it verifies the {\bf Jacobi condition equation} \eqref{Jacobi condition}.
\item[(ii)] A smooth function $C(\vec{M},\vec{\gamma})$ is a Casimir function for the Poisson bracket $\{\cdot,\cdot\}_{\vec{\mu}}$ if and only if $C$ is a function that depends only on the four variables $\gamma_1,\gamma_2,\gamma_3,\vec{M}\cdot \vec{\gamma}$ and it verifies the {\bf Casimir condition equation} \eqref{Casimir condition}.
\item[(iii)] If $\Pi_{\vec{\mu}}$ generates a Poisson bracket, then the system \eqref{sistem general} has the Hamilton-Poisson formulation $(\mathbb{R}^6,\{\cdot,\cdot\}_{\vec{\mu}},H)$, where $H$ is the Hamiltonian function \eqref{Hamiltonian}.
\end{itemize}
\end{thm}

For a gyroscopic torque that verifies the conditions of the above theorem we provide an algorithm that allows us to recover a third constant of motion that we have lost when we added the gyroscopic torque to the generalized Euler-Poisson system \eqref{redus}. This third constant of motion will be a Casimir function for $\{\cdot,\cdot\}_{\vec{\mu}}$, independent of the Casimir $C_1$.

\section{Applications}

We will discuss several well known gyroscopic torques that appear in the literature. We will study these examples in the theoretical frame that we have developed in the previous section.

\subsection{The case $\vec{\mu}=\vec{\mu}(\vec{\gamma})$}

In this particular case when $\vec{\mu}$ depends only on the Poisson variables $\gamma_1,\gamma_2,\gamma_3$ the {\bf Jacobi condition equation} \eqref{Jacobi condition} becomes
\begin{equation}
\label{Jacobi redus}
\vec{\gamma}\cdot \hbox{curl}_{\vec{\gamma}}\vec{\mu}=0,
\end{equation}
which is equivalent with the vector field $\vec{\gamma}\times \vec{\mu}(\vec{\gamma})$ being a solenoidal field, i.e.
$$\hbox{div}_{\vec{\gamma}}(\vec{\gamma}\times \vec{\mu})=0.$$

A large class of examples corresponds to the following solution for equation \eqref{Jacobi redus},
\begin{equation}
\label{mu}
\vec{\mu}(\vec{\gamma})=\psi(\vec{\gamma})\vec{\gamma}+\nabla_{\vec{\gamma}}\varphi(\vec{\gamma}),
\end{equation}
where $\varphi,\psi$ are smooth functions depending only on variables $\gamma_1,\gamma_2,\gamma_3$.\\
For the above $\vec{\mu}$ the {\bf Casimir condition equation} \eqref{Casimir condition} becomes
$$\vec{\gamma}\times ((\partial_4 C)(\psi(\vec{\gamma})\vec{\gamma}+\nabla_{\vec{\gamma}}\varphi(\vec{\gamma}))-\nabla_{\vec{\gamma}}C)=\vec{0}.$$
We search a solution under the hypothesis $\partial_4C=1$. Consequently, we have the following {\bf Casimir condition equation},
$$\vec{\gamma}\times \nabla_{\vec{\gamma}}(\varphi -C)=\vec{0}.$$
A particular solution is given by
\begin{equation}
\label{Casimir particular}
C(\vec{\gamma},\vec{M}\cdot \vec{\gamma})=\vec{M}\cdot \vec{\gamma}+\varphi(\vec{\gamma}).
\end{equation}
This Casimir is functionally independent of the Casimir $C_1=\displaystyle \frac{1}{2}\vec{\gamma}^2$.

In Yehia \cite{Yeh3} it is considered a torque of the form
$$\vec{\mu}(\vec{\gamma})=-(\hbox{div}_{\vec{\gamma}}\vec{l})\vec{\gamma}+\nabla_{\vec{\gamma}}(\vec{l}\cdot \vec{\gamma}),$$
where $\vec{l}=\vec{l}(\vec{\gamma})$ is a smooth function. Choosing $\psi(\vec{\gamma})=-\hbox{div}_{\vec{\gamma}}\vec{l}$ and $\varphi(\vec{\gamma})=\vec{l}\cdot \vec{\gamma}$, $\vec{\mu}$ can be put in the form \eqref{mu} and consequently these torque cases admit a Hamilton-Poisson formulation. The second Casimir is given by $C(\vec{M},\vec{\gamma})=(\vec{M}+\vec{l})\cdot\vec{\gamma}$. This Casimir appears in \cite{Yeh3} as the {\it cyclic integral of motion}.

\medskip

\noindent {\bf The case of gyrostatic torque}

An important case is the gyrostatic torque $\vec{\mathcal{M}}_g$, where $\vec{\mu}(\vec{\gamma})=\vec{\mu}_0$ is a constant vector. A gyrostatic torque can be produced by mechanical rotors or it is the consequence of the liquid action which is contained in rigid
body cavities. An interesting situation for which we have a gyrostatic torque is presented in the Zhukovsky theorem (see \cite{Rum}, pp.
51), where a fluid mass with an initial velocity in a multiply connected cavity also performs an action that is similar to the action of some
rotor attached to the rigid body (the fluid completely fills the cavities of the rigid body and the motion of the fluid is potential). Various
mathematical aspects, for particular functions $U$, are studied in some papers. For a linear function $U$, Gavrilov investigates, in
\cite{Gav}, the integrability of the system. 

A gyrostatic torque can be put in the form \eqref{mu} by choosing $\psi(\vec{\gamma})=0$ and $\varphi(\vec{\gamma})=\vec{\mu}_0\cdot \vec{\gamma}$. From \eqref{Casimir particular} we obtain the second Casimir
$$C(\vec{M},\vec{\gamma})=(\vec{M}+\vec{\mu}_0)\cdot\vec{\gamma}.~~~\blacktriangle$$

\medskip

\noindent {\bf The case of affine gyroscopic torque}

Grioli studies in \cite{Grio} the case of the rotation of a charged rigid body, with a fixed point, situated in a constant
force field. We present some details of the calculation of the torque generated by this force field. We denote by $O$ the fixed point, $q$ the
electric charge, $\vec{\gamma}$ the magnetic field and $D$ the domain of the rigid body. Let $P$ a particle of the rigid body and $\vec{v}(P)$
its velocity. The particle is acted by the Lorenz force
$$\vec{F}=q(\vec{v}\times\vec{\gamma})=q((\vec{\omega}\times\vec{OP})\times\vec{\gamma}).$$
The torque generated by the magnetic field is
\begin{equation*}
\vec{\mathcal{M}}_O=\int_{D}q(P)\vec{OP}\times \vec{F}dv=\vec{\omega}\times\Sigma\vec{\gamma},
\end{equation*}
where the tensor $\Sigma$ is defined by
\begin{equation*}
\Sigma=\int_{D}q(P)\vec{OP}\otimes\vec{OP}dv.
\end{equation*}
It is easy to see that the torque generated by the magnetic field is a gyroscopic torque of the form \eqref{torque} with
$\vec{\mu}=\Sigma\vec{\gamma}$. The matrix associated with the tensor $\Sigma$ is a $3\times 3$ constant symmetric matrix with respect to the body frame $\mathcal{R}$. In \cite{Ple} it is considered a generalized Grioli problem for which the gyroscopic torque is generated by $\vec{\mu}=\Sigma\vec{\gamma}+\vec{e}$, where $\vec{e}$ is a constant vector in the body frame $\mathcal{R}$.

In our case we consider $\vec{\mu}(\vec{\gamma})=A\vec{\gamma}+\vec{\mu}_0$, where $A$ is a $3\times 3$ constant matrix and $\vec{\mu}_0$ is a constant vector, both computed in the body frame.
It is easy to verify that such a torque  $\vec{\mu}(\vec{\gamma})$ satisfies the {\bf Jacobi condition equation} if and only if the matrix $A$ is symmetric. In this case $\vec{\mu}$ can be but in the form \eqref{mu} by choosing $\psi(\vec{\gamma})=0$ and $\varphi(\vec{\gamma})=\displaystyle\frac{1}{2}\vec{\gamma}\cdot A\vec{\gamma} +\vec{\mu}_0\cdot \vec{\gamma}$.

We obtain the second Casimir $C(\vec{M},\vec{\gamma})=\displaystyle\frac{1}{2}\vec{\gamma}\cdot A\vec{\gamma}+(\vec{M}+\vec{\mu}_0)\cdot \vec{\gamma}$. Consequently, the generalized Grioli problem admits a Hamilton-Poisson formulation. $\blacktriangle$

\medskip

\noindent {\bf Nonlinear cases}

In what follows we will analyze several particular cases of systems of type \eqref{sistem general} which are presented in Yehia \cite{Yeh3} and \cite{Yeh4}.\\\\
{\bf a)} We consider a rigid body with the Kovalevskaya configuration ($A=B=2C$), where the potential is given by
$$
U=C(a_1\gamma_1+a_2\gamma_2)-Ck\gamma_3(n+n_1\gamma_1+n_2\gamma_2)-
\frac{1}{2}C(n+n_1\gamma_1+n_2\gamma_2)^2(2\gamma_1^2+2\gamma_2^2+\gamma_3^2)
$$
and the gyroscopic torque is generated by
\begin{align*}
\mu_1&=C(-n\gamma_1-n_1\gamma_1^2+2n_1\gamma_2^2+n_1\gamma_3^2-3n_2\gamma_1\gamma_2)\\
\mu_2&=C(-n\gamma_2+2n_2\gamma_1^2-n_2\gamma_2^2+n_2\gamma_3^2-3n_1\gamma_1\gamma_2)\\
\mu_3&=C(k-3n\gamma_3-5n_1\gamma_1\gamma_3-5n_2\gamma_2\gamma_3),
\end{align*}
where $a_1,a_2,k,n,n_1,n_2$ are real constants.\\
This gyroscopic torque can be put in the form \eqref{mu} with 
$$\varphi=C[k\gamma_3+(n+n_1\gamma_1+n_2\gamma_2)(2\gamma_1^2+2\gamma_2^2+\gamma_3^2)]$$
and 
$$\psi=-C(5n+7n_1\gamma_1+7n_2\gamma_2).$$
Consequently, this case has a Hamilton-Poisson formulation and the second Casimir is given by
$$C(\vec{M},\vec{\gamma})=\vec{M}\cdot\vec{\gamma}+C[k\gamma_3+(n+n_1\gamma_1+n_2\gamma_2)(2\gamma_1^2+2\gamma_2^2+\gamma_3^2)],$$
which coincides with the cyclic integral of motion $I_2$ from \cite{Yeh3}, pp. 342, subsection 4.1.1.
$\diamondsuit$\\\\
{\bf b)} We consider a rigid body with the Kovalevskaya configuration ($A=B=2C$), where the potential is given by
$$
U=C(a_1\gamma_1+a_2\gamma_2)+\frac{\varepsilon}{\sqrt{\gamma_1^2+\gamma_2^2}}
-\frac{1}{2}C\left(n+n_1\gamma_1+n_2\gamma_2+\frac{N}{\sqrt{\gamma_1^2+\gamma_2^2}}\right)^2(2\gamma_1^2+2\gamma_2^2+\gamma_3^2)
$$
and the gyroscopic torque is generated by
\begin{align*}
\mu_1&=C\left(-n\gamma_1-n_1\gamma_1^2+2n_1\gamma_2^2+n_1\gamma_3^2-3n_2\gamma_1\gamma_2+\frac{N\gamma_1}{(\gamma_1^2+\gamma_2^2)^{\frac{3}{2}}}\right)\\
\mu_2&=C\left(-n\gamma_2+2n_2\gamma_1^2-n_2\gamma_2^2+n_2\gamma_3^2-3n_1\gamma_1\gamma_2+\frac{N\gamma_2}{(\gamma_1^2+\gamma_2^2)^{\frac{3}{2}}}\right)\\
\mu_3&=-C\gamma_3\left(3n+5n_1\gamma_1+5n_2\gamma_2+\frac{N\gamma_3}{\sqrt{\gamma_1^2+\gamma_2^2}}\right),
\end{align*}
where $a_1,a_2,\varepsilon,N,n,n_1,n_2$ are real constants.

In Yehia \cite{Yeh3}, pp. 343, subsection 4.1.2, the function
$$I_2=\vec{M}\cdot \vec{\gamma}+C\left(n+n_1\gamma_1+n_2\gamma_2+\frac{N}{\sqrt{\gamma_1^2+\gamma_2^2}}\right)(2\gamma_1^2+2\gamma_2^2+\gamma_3^2)$$
is asserted to be an integral of motion. By direct computation it can be verified that this statement does not hold. We will present a method to solve this incompatibility. 

If we want to keep the gyroscopic torque as above, then  we observe that it verifies the {\bf Jacobi condition equation} \eqref{Jacobi redus} and from the {\bf Casimir condition equation} \eqref{Casimir condition} we obtain
\begin{align*}
C(\vec{M},\vec{\gamma})&=\vec{M}\cdot \vec{\gamma}+C\left(n+n_1\gamma_1+n_2\gamma_2\right)(2\gamma_1^2+2\gamma_2^2+\gamma_3^2)+\\
&+\frac{1}{2}CN\left( \frac{\gamma_3(\gamma_1^2+\gamma_2^2)-2}{\sqrt{\gamma_1^2+\gamma_2^2}}-(\gamma_1^2+\gamma_2^2+\gamma_3^2)\arctan\frac{\gamma_3}{\sqrt{\gamma_1^2+\gamma_2^2}}\right).
\end{align*}

If we want to keep $I_2$ as a constant of motion we will find a new gyroscopic torque generated by $\vec{\mu}'=(\mu_1,\mu_2,\mu_3')$, which can be put in the form \eqref{mu} and for which $I_2$ becomes a Casimir function for the Poisson bracket generated by $\vec{\mu}'$.
We choose
$$\varphi(\vec{\gamma})=C\left(n+n_1\gamma_1+n_2\gamma_2+\frac{N}{\sqrt{\gamma_1^2+\gamma_2^2}}\right)(2\gamma_1^2+2\gamma_2^2+\gamma_3^2)$$
and 
$$\psi(\vec{\gamma})=-C\left(5n+7n_1\gamma_1+7n_2\gamma_2+N\cdot \frac{2(\gamma_1^2+\gamma_2^2)-\gamma_3^2-1}{(\gamma_1^2+\gamma_2^2)^{\frac{3}{2}}}\right).$$
We observe that $\mu_1=\psi(\vec{\gamma})\gamma_1+\frac{\partial \varphi(\vec{\gamma})}{\partial \gamma_1}$, $\mu_2=\psi(\vec{\gamma})\gamma_2+\frac{\partial \varphi(\vec{\gamma})}{\partial \gamma_2}$ and consequently $\mu_3'$ is given by
\begin{align*}
\mu_3'&=\psi(\vec{\gamma})\gamma_3+\frac{\partial \varphi(\vec{\gamma})}{\partial \gamma_3}=\mu_3+CN\gamma_3\frac{\gamma_3(\gamma_1^2+\gamma_2^2)+\gamma_3^2+1}{(\gamma_1^2+\gamma_2^2)^{\frac{3}{2}}}\\
&=-C\gamma_3\left(3n+5n_1\gamma_1+5n_2\gamma_2-N\frac{\gamma_3^2+1}{(\gamma_1^2+\gamma_2^2)^{\frac{3}{2}}}\right).~~~\diamondsuit
\end{align*}

The first, third and fourth cases from Yehia \cite{Yeh4} can be analyzed as in the case {\bf a)}; it is also the situation for the generalization of Lagrange case from Yehia \cite{Yeh4}. The third case from Yehia \cite{Yeh3}, the second and fifth cases from Yehia \cite{Yeh4} can be analyzed as in the case {\bf b)}. $\blacktriangle$

\subsection{The case $\vec{\mu}=\vec{\mu}(\vec{\gamma},\vec{M}\cdot\vec{\gamma})$}

In this section we will study two general subcases of gyroscopic torques  that depend also on $s=\vec{M}\cdot\vec{\gamma}$. We will find a Hamilton-Poisson formulation for these subcases.

\medskip

\noindent {\bf The subcase $\vec{\mu}=a(s)\nabla_{\vec{\gamma}}\varphi (\vec{\gamma})+b(\vec{\gamma},s)\vec{\gamma}$}

For this subcase we consider $a:I\subseteq \mathbb{R}\to\mathbb{R}$, $b:D\times I\subseteq \mathbb{R}^3\times \mathbb{R}\to \mathbb{R}$ and $\varphi :D\subseteq \mathbb{R}^3\to\mathbb{R}$ smooth functions, where $I$ and $D$ are domains in the corresponding spaces.

The gyroscopic torque generated by $\vec{\mu}$ verifies the {\bf Jacobi condition equation}. Indeed, $\vec{\gamma}\cdot \hbox{curl}_{\vec{\gamma}}\vec{\mu}=0$ and
$$
\vec{\mu}\cdot ( \vec{\gamma}\times\partial_4\vec{\mu})=-\vec{\gamma}\cdot ( \vec{\mu}\times \partial_4 \vec{\mu})
=\vec{\gamma}\cdot \left( a(s)\frac{\partial b}{\partial s}(\vec{\gamma},s)\nabla_{\vec{\gamma}}\varphi(\vec{\gamma})\times \vec{\gamma}+
b(\vec{\gamma},s)\frac{\partial a}{\partial s}\vec{\gamma}\times \nabla_{\vec{\gamma}}\varphi(\vec{\gamma})\right)=0.
$$

The {\bf Casimir condition equation} becomes
$$\vec{\gamma}\times \left( (\partial_4C)a(s)\nabla_{\vec{\gamma}}\varphi(\vec{\gamma})-\nabla_{\vec{\gamma}}C\right)=0.$$
We will search for Casimir functions of the form
$$C(\vec{\gamma},s)=D(\vec{\gamma})E(s).$$
Introducing this expression for the Casimir in the above equation we obtain,
$$\vec{\gamma}\times \left(E'(s)a(s)D(\vec{\gamma})\nabla_{\vec{\gamma}}\varphi(\vec{\gamma})-E(s)\nabla_{\vec{\gamma}}D(\vec{\gamma})\right)=0.$$
Sufficient conditions for the above equation to be satisfied are that the functions $E$ and $D$ verify the following equations:
$$a(s)E'(s)=E(s)$$
and
$$D(\vec{\gamma})\nabla_{\vec{\gamma}}\varphi(\vec{\gamma})=\nabla_{\vec{\gamma}}D(\vec{\gamma}).$$
Consequently, we obtain a particular solution 
$$E(s)=\exp\left(\int \frac{1}{a(s)}ds\right)$$
and
$$D(\vec{\gamma})=\exp \varphi(\vec{\gamma}).$$
Thus, a Casimir function for $\Pi_{\vec{\mu}}$ is given by
\begin{equation}
\label{Casimir case 1}
C(\vec{\gamma},s)=\int \frac{1}{a(s)} ds+\varphi(\vec{\gamma}).
\end{equation}

This case is a generalization of the gyroscopic torque given by relation \eqref{mu} for $a(s)=1$ and $b(\vec{\gamma},s)=\psi(\vec{\gamma})$. The Casimir \eqref{Casimir particular} is a particular case of the Casimir function \eqref{Casimir case 1}. $\blacktriangle$

\medskip

\noindent {\bf Gyroscopic torques along one axis of inertia}

Without losing the generality we will study gyroscopic torques generated by $\vec{\mu}=(0,0,\mu_3(\vec{\gamma},s))$.

The {\bf Jacobi condition equation} for $\vec{\mu}$ of this form reduces to
\begin{equation}
\label{Jacobi case 2}
\gamma_1\frac{\partial \mu_3}{\partial \gamma_2}-\gamma_2\frac{\partial \mu_3}{\partial \gamma_1}=0,
\end{equation}
which has the solution $\mu_3=\mu_3(r,\gamma_3,s)$, where $r:=\displaystyle\frac{1}{2}(\gamma_1^2+\gamma_2^2)$.

The {\bf Casimir condition equation} is equivalent with the following system of equations:
$$\left\{\begin{array}{l}
\gamma_2\frac{\partial C}{\partial\gamma_1}-\gamma_1\frac{\partial C}{\partial \gamma_2}=0\\
\gamma_2\left(\frac{\partial C}{\partial s}\mu_3-\frac{\partial C}{\partial\gamma_3}\right)+\gamma_3\frac{\partial C}{\partial\gamma_2}=0\\
\gamma_1\left(\frac{\partial C}{\partial s}\mu_3-\frac{\partial C}{\partial\gamma_3}\right)+\gamma_3\frac{\partial C}{\partial\gamma_1}=0.
\end{array}\right.$$
From the first equation we obtain that the Casimir function has the form 
$$C=C(r,\gamma_3,s).$$
Substituting this form in the last two equations we obtain
\begin{equation}
\label{Casimir case 2}
\mu_3\frac{\partial C}{\partial s}-\frac{\partial C}{\partial \gamma_3}+\gamma_3\frac{\partial C}{\partial r}=0.
\end{equation}
For $\mu_3$ of the particular form 
$$\mu_3(r,\gamma_3,s)=\beta(\gamma_3)\delta(s)$$
we have a Casimir function given by
$$C(r,\gamma_3,s)=\int \frac{1}{\delta(s)}ds+\int \beta(\gamma_3)d\gamma_3$$
as a solution of \eqref{Casimir case 2}.

In Borisov and Mamaev \cite{Bor} it is introduced a gyroscopic torque generated by $\mu_1=\mu_2=0$ and $\mu_3=\displaystyle\frac{\vec{M}\cdot\vec{\gamma}}{\gamma_3}$ for the Kovalevskaya configuration of the rigid body and the potential $U=\alpha (\gamma_1^2-\gamma_2^2)$, where $\alpha$ is a constant.
This choice for $\vec{\mu}$ verifies \eqref{Jacobi case 2} and the function $C(r,\gamma_3,\vec{M}\cdot\vec{\gamma})=\gamma_3 s=\gamma_3(\vec{M}\cdot \vec{\gamma})$ verifies \eqref{Casimir case 2}
and thus is a Casimir function for $\Pi_{\vec{\mu}}$. $\blacktriangle$


\begin{thebibliography}{99}

\bibitem{Yeh3} H.M. Yehia, {\it Regular and Chaotic Dynamics} {\bf 3} (2003) 337-348.
\bibitem{Yeh4} H.M. Yehia, {\it J. Phys. A: Math. Gen} {\bf 30} (1997) 7269-7275.
\bibitem{Bor} A.V. Borisov, I.S. Mamaev, {\it Regular and Chaotic Dynamics} {\bf 3} (1997) 72-89.
\bibitem{Chang1} D.E. Chang, J.E. Marsden, {\it SIAM J. Control Optim.} {\bf 1} (2004) 277-300.
\bibitem{Chang2} D.E. Chang, A.M. Bloch, N.E. Leonard, J.E. Marsden, C.A. Woolsey, {\it ESAIM. Control, Optimisation and Calculus of Variations} {\bf 8} (2002) 393-422.
\bibitem{Rum} N.N. Moiseyev, V.V. Rumyantsev, {\it Dynamic Stability of Bodies Containing Fluid}, Springer-Verlag, Berlin-Heidelberg-New
York, 1968.
\bibitem{Gav} L. Gavrilov, {\it Compositio Mathematica} {\bf 3} (1992) 257-291.
\bibitem{Grio} G. Grioli, {\it Rendiconti del Seminario Matematico della Universita di Padova} {\bf 27} (1957) 90-102.
\bibitem{Ple} Y.D. Pleshakov, {\it Doklady Physics} {\bf 4} (2007) 225-227.





\end{thebibliography}
\end{document}